\documentclass[11pt,draft]{article}
\usepackage{amsmath}
\usepackage{amssymb}
\usepackage{enumerate}
\usepackage{color}
\usepackage{xfrac}

\newtheorem{thm}{Theorem}[section]
\newtheorem{defn}[thm]{Definition}
\newtheorem{lemma}[thm]{Lemma}
\newtheorem{prop}[thm]{Proposition}
\newtheorem{rem}[thm]{Remark}
\newtheorem{cor}[thm]{Corollary}
\newtheorem{ex}[thm]{Example}

\DeclareMathOperator{\Log}{Log}
\DeclareMathOperator{\Arg}{Arg}

\begin{document}

\title{Nielsen's beta function and some infinitely divisible distributions}

\author{Christian Berg, Stamatis Koumandos and Henrik L. Pedersen\footnote{Research supported by grant DFF–4181-00502 from The Danish Council for Independent Research $|$ Natural Sciences}
}

\date{\today}
\maketitle

\begin{abstract}
We show that a large collection of special functions, in particular Nielsen's beta function, are generalized Stieltjes functions of order 2, and therefore logarithmically 
completely monotonic. This includes the Laplace transform of functions of the form $xf(x)$, where $f$ is itself the Laplace transform of a sum of dilations and translations
of periodic functions. Our methods are also applied to ratios of Gamma functions, and to the remainders in asymptotic expansions of the double Gamma function of Barnes. 
\end{abstract}


\section{Introduction}
Nielsen's beta function is the meromorphic function defined by 
$$
\beta(z)=\sum_{n=0}^{\infty}\frac{(-1)^n}{z+n},
$$
cf. \cite[p.~101]{N} and \cite[p.~907]{Gradshteyn}. In the right half plane it admits the integral representation 
\begin{equation}
 \label{eq:beta-int}
\beta(z)=\int_0^{\infty}\frac{e^{-zt}}{1+e^{-t}}\, dt,
\end{equation}
from which it immediately follows that $\beta$ is a completely monotonic function on $(0,\infty)$. 
It is a consequence of our results that $\beta$ is also a so-called logarithmically completely monotonic function.
For the reader's convenience we give some background material on the relation between completely monotonic
functions, logarithmically completely monotonic functions and generalized Stieltjes functions.

A completely monotonic function $f$ is a $C^{\infty}$-function on $(0,\infty)$ such that $(-1)^nf^{(n)}(x)\geq 0$ for all $n\geq 0$ and all $x>0$. 
These functions are characterized in Bernstein's theorem as Laplace transforms of positive measures: $f$ is completely monotonic if and only if there exists 
a positive measure $\mu$ on $[0,\infty)$ such that $t\mapsto e^{-xt}$ is integrable w.r.t.\ $\mu$ for all $x>0$ and
$$
f(x)=L(\mu)(x)=\int_0^{\infty}e^{-xt}\, d\mu(t),
$$
cf.\ \cite[p.~161]{widder}.

For a given  $\lambda >0$ a function $f:(0,\infty)\to \mathbb R$ is called a generalized Stieltjes function of order $\lambda$ if 
\begin{equation}
\label{eq:def-S}
f(x)=\int_0^{\infty}\frac{d\mu(t)}{(x+t)^{\lambda}}+c,
\end{equation}
where $\mu$ is a positive measure on $[0,\infty)$ making the integral converge for $x>0$ and $c\geq 0$. The set of generalized Stieltjes functions of order $\lambda$ is denoted $\mathcal S_\lambda$. 
For additional information on these classes see e.g.\ \cite{kp2}.

The class of ordinary Stieltjes functions is the class of generalized Stieltjes functions of order $1$.

We remark that $f$ is a generalized Stieltjes function of order $\lambda$ of the form \eqref{eq:def-S} if and only if 
\begin{equation}
\label{eq:sokalabsolutely}
f(x)=\frac{1}{\Gamma(\lambda)}\int_0^{\infty}e^{-xt}t^{\lambda-1}\kappa(t)\, dt+c, \quad x>0,
\end{equation}
where $\kappa$ is a completely monotonic function. In the affirmative case, 
$
\kappa(t)=\int_0^{\infty}e^{-ts}\, d\mu(s).
$
See \cite[Lemma 2.1]{kp2}. This characterization shows also that
$\mathcal S_{\lambda_1}\subseteq \mathcal S_{\lambda_2}$ for $\lambda_1<\lambda_2$.

Section \ref{sec:proof} contains the proof of Theorem \ref{thm:main} below, covering Nielsen's beta function, based on a result of Kristiansen and Steutel.
Some additional material on the relation between generalized Stieltjes functions and so-called logarithmically completely monotonic functions (defined below) is also given.

In Section \ref{sec:periodic} we describe a method in which we start with sums $\varphi$ of certain translations and dilations of an even, non-negative and periodic function and then consider 
the Laplace transform of
the function $\sigma(x)=xL(\varphi)(x)$. From \eqref{eq:sokalabsolutely} we thus obtain that $L(\sigma)$ is in $\mathcal S_2$ and has the same representing measure in 
\eqref{eq:def-S} as $\varphi$ has as a completely monotonic function. As corollaries we obtain analogues of results of Kosaki (see \cite{kosaki}) on 
infinitely divisible functions in relation to the Fourier transform. 

Section \ref{sec:gamma} deals with ratios of Gamma functions and contains an improvement of results of Bustoz and Ismail (see \cite{bustoz-ismail}). In Section \ref{sec:barnes} we improve earlier results for 
the remainders in asymptotic expansions of the double Gamma function of Barnes.

Finally we consider in Section \ref{sec:cesaro} generalized Stieltjes functions related to in particular meromorphic functions of special type having poles at the non-positive integers.

A $C^{\infty}$-function $f:(0,\infty)\to (0,\infty)$ is said to be logarithmically completely monotonic if $-(\log f(x))'=-f'(x)/f(x)$ is completely monotonic. 
It is easy to see that a logarithmically completely monotonic function is also completely monotonic. The class
of logarithmically completely monotonic functions is denoted by $\mathcal L$ and was characterized by Horn (\cite{horn}) through Theorem \ref{thm:horn}. See also \cite{berg1,berg2}.
Because of (ii) the functions of $\mathcal L$ were called infinitely divisible completely monotonic functions by Horn.
\begin{thm}
\label{thm:horn}
 The following conditions are equivalent for a function $f: (0,\infty)\to (0,\infty)$:
 \begin{enumerate}
  \item[(i)] $f\in \mathcal L$,
  \item[(ii)] $f^{\alpha}$ is completely monotonic for all $\alpha >0$,
  \item[(iii)] $f^{1/n}$ is completely monotonic for all $n=1,2,\ldots$
 \end{enumerate}
\end{thm}
In \cite[p.~90]{steutel-van-harn} it is noticed that a positive differentiable function $\pi$ on $(0,\infty)$ with $\lim_{x\to 0+}\pi(x)=1$ is an infinitely divisible probability Laplace-Stieltjes transform
if and only if its so-called $\rho$-function (which is the negative of its logarithmic derivative) is completely monotonic. However, the general concept of logarithmic complete monotonicity is not 
mentioned explicitly in \cite{steutel-van-harn}. 

\begin{defn}
A measure or function on the positive half line is called infinitely divisible if its Laplace transform belongs to the class $\mathcal L$.
\end{defn}
We shall first investigate functions of the form 
\begin{equation}\label{eq:def-f}
f(z)=\sum_{n=0}^{\infty}\frac{(-1)^n}{(z+a_n)^{\lambda}},
\end{equation}
where $0\leq a_0\leq a_1\leq\cdots \to \infty$, and where $\lambda>0$. Notice that the series converges uniformly over compact subsets of the cut plane $\mathbb C\setminus (-\infty,0]$. 
Our first main result is the following.
\begin{thm}
\label{thm:main}
 Let $f$ be of the form \eqref{eq:def-f} for $\lambda\in (0,1]$. Then $f$ is a generalized Stieltjes function of order $\lambda+1$ with $c=0$ and representing measure
 $$
 \mu=\lambda\sum_{n=0}^{\infty}dt|_{(a_{2n}, a_{2n+1})}.
 $$
 In particular $f$ is a logarithmically completely monotonic function. Furthermore, it admits the representation
 $$
 f(x)=\frac{1}{\Gamma(\lambda+1)}\int_0^{\infty}e^{-xt}t^{\lambda}\varphi(t)\, dt, \quad x>0,
 $$
 where $\varphi(t)=\sum_{n=0}^{\infty}(e^{-a_{2n}t}-e^{-a_{2n+1}t})/t$ is completely monotonic.
\end{thm}
Nielsen's beta function is a special case of \eqref{eq:def-f}. By Theorem \ref{thm:main} we then get:
\begin{cor}
\label{cor:m_c}
There exists a family of measures $\{m_c\}_{c>0}$ on $[0,\infty)$ such that 
\begin{enumerate}
 \item[(a)] $\beta(z)^c=\int_0^\infty e^{-zt}\,dm_c(t),\quad \Re z>0$,
 \item[(b)] $m_c*m_d=m_{c+d}$,
 \item[(c)] $dm_1(t)=dt/(1+e^{-t})$.
\end{enumerate}
\end{cor}
Let $a>0$ and define
$$
\nu_a(t)=\frac{e^{-at}}{\beta(a) (1+e^{-t})}=\frac{1}{2\beta(a)}\frac{e^{-t(a-1/2)}}{\cosh(t/2)}
$$
It follows from Theorem \ref{thm:main} that any of these probability densities is infinitely divisible.  
For $a=1/2,1,3/2$ we find
$$
\nu_{1/2}(t)=\frac{1}{\pi \cosh(t/2)}, \ \nu_1(t)=\frac{1}{(\log 2)(e^t+1)}, \ \nu_{3/2}(t)=
\frac{e^{-t}}{(4-\pi)\cosh(t/2)}.
$$
We could not find any reference indicating that these probability densities are infinitely divisible.

The integral representation \eqref{eq:beta-int} of $\beta$ may be written as
\begin{equation}
 \label{eq:beta-cosh}
\beta\left(\frac{x+1}{2}\right)=
\int_{0}^{\infty}e^{-xt}\,\frac{1}{\cosh t}\,dt,
\end{equation}
and the logarithmic complete monotonicity of this function corresponds to the infinite divisibility of $\nu_{1/2}$ above.

\section{Generalized Stieltjes functions of order 2}
\label{sec:proof}

\noindent
{\it Proof of Theorem \ref{thm:main}.} For $x>0$ we have 
\begin{align*}
f(x)=\sum_{n=0}^{\infty}\frac{(-1)^n}{(x+a_n)^{\lambda}}&=\sum_{n=0}^{\infty}\frac{1}{(x+a_{2n})^{\lambda}}-\frac{1}{(x+a_{2n+1})^{\lambda}}\\
&=\lambda\sum_{n=0}^{\infty}\int_{a_{2n}}^{a_{2n+1}}\frac{dt}{(x+t)^{\lambda +1}} =\int_{0}^{\infty}\frac{d\mu(t)}{(x+t)^{\lambda +1}},
\end{align*}
where 
$$
\mu=\lambda\sum_{n=0}^{\infty}dt|_{(a_{2n}, a_{2n+1})}.
$$
Therefore, $f$ is a generalized Stieltjes function of order $\lambda +1$, and hence also of order $2$. Such functions 
are logarithmically 
completely monotonic according to Theorem \ref{thm:Kristiansen} below.
The representation as a Laplace transform follows easily from \eqref{eq:sokalabsolutely}. \hfill $\square$

We note that a result similar to Theorem \ref{thm:main} holds for functions given by a finite sum. Indeed, considering for $N\geq 0$
$$
f_N(x)=\sum_{n=0}^{N}\frac{(-1)^n}{(x+a_n)^{\lambda}},
$$
it follows that 
$$
f_{2N+1}(x)=\int_{0}^{\infty}\frac{d\mu_{2N+1}(t)}{(x+t)^{\lambda +1}},\quad \text{with}\ \mu_{2N+1}=\lambda\sum_{n=0}^{N}dt|_{(a_{2n}, a_{2n+1})}
$$
and 
$$
f_{2N}(x)=\int_{0}^{\infty}\frac{d\mu_{2N}(t)}{(x+t)^{\lambda +1}},\quad \text{with}\ \mu_{2N}=\lambda\sum_{n=0}^{N-1}dt|_{(a_{2n}, a_{2n+1})}+\lambda dt|_{(a_{2N}, \infty)}.
$$

\begin{thm}
\label{thm:Kristiansen}
Any generalized Stieltjes function of order $2$ is logarithmically completely monotonic, i.e.\ $\mathcal S_2\subset\mathcal L$.
\end{thm}
This theorem is a deep result. It was conjectured by Steutel in 1970 (see \cite[p.~43]{steutel}) that all functions of the form 
$$
\int_0^{\infty}\frac{t^2}{(x+t)^2}\, d\sigma(t),
$$
where $\sigma$ is a probability measure on $(0,\infty)$, are infinitely divisible completely monotonic functions, or equivalently that all probability densities of the form $xg(x)$, with $g$ being completely 
monotonic, are infinitely divisible. Furthermore, Steutel (see \cite{steutel1}) showed in 1980 that this conjecture would be verified provided another conjecture about the zero distribution of certain rational functions in the complex plane would hold true.
Kristiansen proved the latter conjecture in 1994 (see \cite{kristiansen}). 

The result by Steutel and Kristiansen is easily extended to all generalized Stieltjes functions of order $2$ by a limit argument: If $f\in \mathcal S_{\lambda}$ has the representation 
$$
f(x)=\int_{0}^{\infty}\frac{d\mu(t)}{(x+t)^{2}}+c
$$
then consider $f_n(x)=\int_{0}^{\infty}d\mu_n(t)/(x+t)^{2}$ with $\mu_n=\mu(\{0\})\epsilon_{1/n} + \mu|_{(0,n)} + cn^2\epsilon_n$ ($\epsilon_a$ denoting the point mass at $a$)
 and notice that by Steutel and Kristiansens result $f_n$ is in $\mathcal L$.  Since $\mathcal L$ is closed under pointwise 
convergence, $f$ is also in $\mathcal L$. In this way Theorem \ref{thm:Kristiansen} is proved.

\begin{rem}
The function $f(x)=1/\sinh(x)$ is logarithmically completely monotonic. Indeed, a direct computation shows that 
$$
-f'(x)/f(x)=L\left(\epsilon_0+2\sum_{k=1}^{\infty}\epsilon_{2k}\right)(x).
$$
However, $f$ is not in any of the classes $\mathcal S_{\lambda}$ since $x^{\lambda}f(x)\to 0$ as $x\to \infty$. 

Furthermore $f=L(\nu)$, with $\nu=2\sum_{k=0}^\infty \epsilon_{2k+1}$. 
An old formula of Hamburger (see \cite{Hamburger} and \cite[p.~486]{Gradshteyn}), reads (in a modified form)
$$
\left(x\prod_{k=1}^n \left(1+x^2/(k\pi)^2\right)\right)^{-1}=\frac{2^n}{\binom{2n}{n}}\int_0^\infty e^{-xt}(1-\cos(t\pi))^n\ dt.
$$
The left-hand side is the partial product of $f$ and the measure on the right hand side converges vaguely to the discrete measure $\nu$.
\end{rem}

The subclass of those functions $g$ in $S_{\lambda}$ for which $\lim_{x\to \infty}g(x)=0$ is denoted by $T_{\lambda}$ and it follows immediately that 
$S_{\lambda}= T_{\lambda}\oplus [0,\infty)$. The example in \cite[p.~409]{steutel-van-harn} shows that for any $r>2$ the set $S_r\setminus \mathcal L$ is non empty. 
The following example shows that one can even find concrete functions $g$ in $T_r\setminus \mathcal L$.

\begin{ex}
\label{ex:r>2}
Let $r>2$ and consider, for $c>0$,
$$
g_c(z)=\frac{c}{1+c}\frac{1}{z^r}+\frac{1}{1+c}\frac{1}{(z+1)^r},\quad \Re z>0.
$$
(Here, the power is defined in terms of the principal logarithm.) We claim that for some $c>0$ the function $g_c$ is not in $\mathcal L$. The conformal mapping $\varphi(z)=z/(z+1)$
maps the right half plane onto the open disk of radius $1/2$ centered at $1/2$ and therefore $z\mapsto \varphi(z)^r$ maps the right half plane onto a region covering a non empty
interval $I_r$ of the negative real axis (at least twice). Choosing now $c>0$ with $-c\in I_r$ there must exist $z_c\in \{ z\, | \, \Re z>0\}$ such that
$\varphi(z_c)^r=-c$. Since $\varphi(z)^r=z^r/(z+1)^r$ it therefore follows that $g_c(z_c)=0$. This gives a contradiction since a logarithmically completely monotonic function
cannot have any zeros in the right half plane. 
\end{ex}
Note that Example \ref{ex:r>2} shows that $\mathcal L$ is not stable under addition.
\begin{rem}
 The subclass $T_{\lambda}$ is dense in $S_{\lambda}$ w.r.t.\ pointwise convergence. Indeed, putting $\mu_n=n^{\lambda}\epsilon_n$, it follows that $ \int_0^{\infty}d\mu_n(t)/(x+t)^{\lambda}\to 1$
 as $n\to \infty$ and this shows the denseness. As a consequence the example in \cite{steutel-van-harn} mentioned above also implies that for any $r>2$ there must exist functions $g$ in $T_r\setminus \mathcal L$. 
\end{rem}
\begin{rem}
 Euler's digamma function $\psi=\Gamma' / \Gamma$ is given by
 $$
 \psi(x)=-\gamma +\sum_{n=0}^{\infty}\left(\frac{1}{n+1}-\frac{1}{n+x}\right).
 $$
We have
$$
\psi'(x)=\sum_{n=0}^\infty\frac{1}{(x+n)^2}=\int_0^\infty\frac{d\mu(t)}{(x+t)^2},
$$
where $\mu=\sum_{n=0}^\infty \epsilon_n$. Clearly $L(\mu)(t)=(1-e^{-t})^{-1}$.
It follows that
$$
\psi'=L\left(\frac{t}{1-e^{-t}}\right)\in\mathcal S_2\subset \mathcal L.
$$
The measure with density on the half line given by 
$$
\frac{t}{1-e^{-t}}=\frac{te^{t/2}}{2\sinh(t/2)}
$$
is thus infinitely divisible. For $a>0$ define
$$
\tau_a(t)=\frac{te^{-at}}{\psi'(a)(1-e^{-t})}.
$$
This is an infinitely divisible probability density on the half line and for $a=1/2,1,3/2$ we find
$$
\tau_{1/2}(t)=\frac{t}{\pi^2\sinh(t/2)},\; \tau_1(t)=\frac{6}{\pi^2}\frac{t}{e^t-1},\;
\tau_{3/2}(t)=\frac{te^{-t}}{(\pi^2-8)\sinh(t/2)}.
$$
\end{rem}

\section{Infinitely divisible functions related to periodic functions}
\label{sec:periodic}
In this section we obtain a number of infinitely divisible functions by studying Laplace transforms of a class of positive and periodic functions. If $\varphi:[0,\infty)\to \mathbb R$ is periodic with
period $T>0$ then it is easy to see that
\begin{equation}
\label{eq:simple-periodic}
L(\varphi)(x)=\frac{1}{1-e^{-xT}}\int_0^{T}\varphi(t)e^{-xt}\, dt.
\end{equation}
We first consider the case of a periodic and piecewise constant function $\varphi$. 
\begin{thm}
\label{thm:discrete}
Let $a_{1},a_{2}, \ldots, a_{n}$ be nonnegative real numbers, $a_{n}=a_{1}$,  $0<\lambda_{1}< \lambda_{2}<\ldots<\lambda_{n}$, $\alpha>0$ and $0\leq \beta\leq \lambda_{1}$. 
\begin{enumerate}
\item[(i)] The 
Laplace transform of the function $\sigma$ given by
$$
\sigma(x)=a_{1}+\cosh\left(\frac{\beta}{\alpha}\,x\right)\,\sum_{k=1}^{n-1}(a_{k}-a_{k+1})\frac{1-e^{-\frac{\lambda_{k}}{\alpha}\,x}}{1-e^{-\frac{\lambda_{n}}{\alpha}\,x}}
$$
is a generalized Stieltjes function of order 2 and hence $\sigma$ is infinitely divisible.
\item[(ii)] The 
Laplace transform of the function $\tau$ given by
$$
\tau(x)=\max_{1\leq j\leq n}\{a_{j}\}\pm2\sinh\left(\frac{\beta}{\alpha}\,x\right)\,\sum_{k=1}^{n-1}(a_{k}-a_{k+1})\frac{1-e^{-\frac{\lambda_{k}}{\alpha}\,x}}{1-e^{-\frac{\lambda_{n}}{\alpha}\,x}}
$$
is a generalized Stieltjes function of order 2 and hence $\tau$ is infinitely divisible.
\end{enumerate}
\end{thm}
{\it Proof.}
Consider the function
$$
\varphi(t)=\sum_{k=1}^{n}a_{k}\chi_{[\lambda_{k-1}, \lambda_{k})}(t), \;\;\text{for all}\;\; t\in[0,T),
$$
where $T=\lambda_{n}$, $\lambda_{0}=0$, and $\chi_{I}$ denotes the indicator function of $I$.
Extend $\varphi$ to $[0,\infty)$ by $T$-periodicity and then to $\mathbb R$ by $\varphi(t)=\varphi(-t)$ so that $\varphi$ is a nonnegative $T$-periodic and even function. Denoting  the 
Laplace transform of $\varphi$ by $F$ we see that
\begin{align*}
F(x)&=\frac{1}{1-e^{-\lambda_{n}x}}\int_{0}^{\lambda_{n}}\varphi(t)e^{-xt}\,dt\\
&= \frac{1}{x}\,\frac{1}{1-e^{-\lambda_{n}x}}\sum_{k=1}^{n}a_{k}(e^{-\lambda_{k-1}x}-e^{-\lambda_{k}x})\nonumber\\[+3pt]
&= \frac{1}{x}\,\left(a_{n}+\sum_{k=1}^{n-1}(a_{k}-a_{k+1})\frac{1-e^{-\lambda_{k}x}}{1-e^{-\lambda_{n}x}}\right).
\end{align*}
Here the last equality is obtained by a summation by parts, using  $a_{n}=a_{1}$.
For $\alpha>0$ and $0\leq \beta\leq \lambda_{1}$ we obtain
\begin{equation}
\int_{0}^{\infty}\varphi(\alpha\,t-\beta)e^{-xt}\,dt=\frac{a_{n}}{x}-\frac{a_{n}}{x}\,e^{-\frac{\beta}{\alpha}x}+\frac{1}{\alpha}F\left(\frac{x}{\alpha}\right)\,e^{-\frac{\beta}{\alpha}x}.
\end{equation}
Since $a_{n}=a_{1}$ we similarly have
\begin{equation}
\int_{0}^{\infty}\varphi(\alpha\,t+\beta)e^{-xt}\,dt=\frac{a_{n}}{x}-\frac{a_{n}}{x}\,e^{\frac{\beta}{\alpha}x}+\frac{1}{\alpha}F\left(\frac{x}{\alpha}\right)\,e^{\frac{\beta}{\alpha}x}.
\end{equation}
Combining these relations it follows that
\begin{align}
&\frac{1}{2}\int_{0}^{\infty}\left(\varphi(\alpha\,t+\beta) +\varphi(\alpha\,t-\beta)\right)\,e^{-xt}\,dt=\nonumber\\[+5pt]
&\frac{1}{x}\left(a_{1}+\cosh\left(\frac{\beta}{\alpha}\,x\right)\, \sum_{k=1}^{n-1}(a_{k}-a_{k+1})\frac{1-e^{-\frac{\lambda_{k}}{\alpha}\,x}}{1-e^{-\frac{\lambda_{n}}{\alpha}\,x}}\right).\label{eq:alpha-beta}
\end{align}
Hence $\sigma(x)/x$ is completely monotonic ensuring that $L(\sigma)$ is in $\mathcal S_2$, and hence in $\mathcal L$. This completes the proof of (i).

To prove (ii) we let $A=\max_{1\leq j\leq n}\{a_{j}\}$ and obtain
\begin{align*}
&\int_{0}^{\infty}\left(A\pm (\varphi(\alpha\,t+\beta) -\varphi(\alpha\,t-\beta))\right)\,e^{-xt}\,dt=\\
&\frac{1}{x}\left(A\pm2\sinh\left(\frac{\beta}{\alpha}\,x\right)\, \sum_{k=1}^{n-1}(a_{k}-a_{k+1})\frac{1-e^{-\frac{\lambda_{k}}{\alpha}\,x}}{1-e^{-\frac{\lambda_{n}}{\alpha}\,x}}\right).
\end{align*}
Since $A\pm (\varphi(\alpha\,t+\beta) -\varphi(\alpha\,t-\beta))$ is nonnegative, $\tau(x)/x$ is also completely monotonic.\hfill $\square$

\begin{cor}
\label{cor:+}
The functions 
$$
\sigma_1(x)=1-\frac{\cosh(\nu x)}{2\cosh(x)}\quad \text{and}\quad  \sigma_2(x)=\frac{\cosh(\nu x)}{\cosh(x)},\;\; \nu\in[0, 1]
$$
are infinitely divisible. Their Laplace transforms belong to $\mathcal S_2$ and are given by
\begin{align*}
L(\sigma_1)(x)&=\frac{1}{x}-\frac{1}{4}\left(\beta\left(\frac{x+1-\nu}{2}\right)+\beta\left(\frac{x+1+\nu}{2}\right)\right),\\
L(\sigma_2)(x)&=\frac{1}{2}\left(\beta\left(\frac{x+1-\nu}{2}\right)+\beta\left(\frac{x+1+\nu}{2}\right)\right).
\end{align*}
\end{cor}
{\it Proof.} 
An application of (i) of Theorem \ref{thm:discrete} with
$n=3$, $\lambda_{1}=\pi/2$, $\lambda_{2}=3\pi/2$, $\lambda_{3}=2\pi$,	 $a_{1}=a_{3}=1$ $a_{2}=0$,
$\alpha=\pi/2$ and $\beta=\pi\nu/2$, corresponding to (the even and periodic extension of) $\varphi(t)=\chi_{[0,\pi/2)}(t)+\chi_{[3\pi/2, 2\pi)}(t)$ for $t\in [0,2\pi)$ and 
$$
F(x)=L(\varphi)(x)=\frac{1}{x}\left(1-\frac{1}{2\cosh(\pi x/2)}\right)
$$
yields infinite divisibility of 
$$
\sigma_1(x)=1-\frac{\cosh(\nu x)}{2\cosh(x)}\,.
$$
Consider next the function
$$
p(t)=2-\varphi\left(\frac{\pi}{2}t-\frac{\pi\nu}{2}\right)-\varphi\left(\frac{\pi}{2}t+\frac{\pi\nu}{2}\right).
$$
It is clear that $p(t)\geq 0$ and by applying \eqref{eq:alpha-beta} we obtain
\begin{equation}
\int_{0}^{\infty} e^{-xt}\,p(t)\,dt=\frac{\cosh(\nu x)}{x\cosh(x)}.
\end{equation}
Therefore the function 
$$
\sigma_2(x)=\frac{\cosh(\nu x)}{\cosh(x)}
$$
is infinitely divisible. The Laplace transforms are computed using \eqref{eq:beta-cosh}.\hfill $\square$
\begin{rem}
 When $\nu=0$ in the above corollary, $\sigma_2(x)=1/\cosh(x)$ and we revisit the density function $\nu_{1/2}$ mentioned right after Corollary  \ref{cor:m_c}.
\end{rem}

\begin{cor}
\label{cor:+-}
The functions 
$$
\tau(x)=1\pm\frac{\sinh(\nu x)}{\cosh(x)}
\quad \nu\in[0, 1]
$$
are infinitely divisible. Their Laplace transforms belong to $\mathcal S_2$ and are given by
$$
L(\tau)(x)=\frac{1}{x}\pm \frac{1}{2}\left(\beta\left(\frac{x+1-\nu}{2}\right)-\beta\left(\frac{x+1+\nu}{2}\right)\right).
$$ 
\end{cor}
{\it Proof.}
Applying (ii) of Theorem \ref{thm:discrete} with
$n=3$, $\lambda_{1}=\pi/2$, $\lambda_{2}=3\pi/2$, $\lambda_{3}=2\pi$,	 $a_{1}=a_{3}=1$ $a_{2}=0$,
$\alpha=\pi/2$ and $\beta=\pi\nu/2$ yields the infinite divisibility.\hfill $\square$

\begin{rem}
These corollaries are analogues of results in \cite{kosaki}. In that paper infinite divisible functions with respect to the Fourier transform are studied via a general
result (\cite[Theorem 2]{kosaki}) about entire functions. That result does not carry over to our situation. 
\end{rem}
\begin{rem}
 \label{rem:psi-beta}
 Returning to Euler's digamma function we notice
 $$
 \beta(x)=\frac{1}{2}\left(\psi\left(\frac{x+1}{2}\right)-\psi\left(\frac{x}{2}\right)\right)
 $$
 and it can be seen that such a difference is in $\mathcal S_2$.
\end{rem}

We continue to the case where $\varphi$ is continuous and piecewise smooth. The following lemma is easily obtained from \eqref{eq:simple-periodic} by partial integration. 
\begin{lemma}\label{lemma:continuous}
Let $\varphi$ be a non-negative,  $T$-periodic  function. Suppose furthermore that $\varphi$ is continuous and piecewise smooth on $[0,T)$. Then
\begin{align*}
L(\varphi)(x)
&= \frac{1}{x}\,\varphi(0)+\frac{1}{x}\,\,\frac{1}{1-e^{-Tx}}\,  \int_{0}^{T}\varphi^{\prime}(t)e^{-xt}\,dt.
\end{align*}
\end{lemma}
The next corollary is also an analogue of a result in \cite{kosaki}. 
\begin{cor}
\label{cor:cauchy}
The function
$$
\omega(x):=\frac{1+bx^2}{1+ax^2}
$$
is infinitely divisible for $a\geq b\geq 0$, its Laplace transform is in $\mathcal S_2$ and is given by 
$$
  L\left(\omega\right)(x)=\frac{b}{a}\frac{1}{x} +\left(1-\frac{b}{a}\right)\,\frac{1}{\sqrt{a}}\,
  \left(\operatorname{ci}(\frac{x}{\sqrt{a}})\,\sin \frac{x}{\sqrt{a}}- \operatorname{si}(\frac{x}{\sqrt{a}})\,\cos \frac{x}{\sqrt{a}}\right),  
$$
where 
\begin{align*}
     \operatorname{si}(x)=-\int_{x}^{\infty}\frac{\sin t}{t}\, dt\quad \text{and} \quad \operatorname{ci}(x)=-\int_{x}^{\infty}\frac{\cos t}{t}\, dt.
\end{align*}
\end{cor}
{\it Proof.}
An application of Lemma \ref{lemma:continuous} for
$$
\varphi(t)=1-\left(1-\frac{b}{a}\right)\cos \frac{t}{\sqrt{a}}
$$
yields 
$$
L(\varphi)(x)=\frac{1+bx^2}{x(1+ax^2)}.
$$
We observe that 
$$
  \omega(t)=\frac{b}{a} +\left(1-\frac{b}{a}\right)\,\frac{1}{at^2+1},
  $$
from which its Laplace transform can easily be calculated (cf.\ \cite[p.~342]{Gradshteyn}).\hfill $\square$

We have in the particular case $b=0$, $a=1$ the Cauchy density function and 
$$
F(x):=L\left(\frac{1}{1+t^2}\right)(x)= \operatorname{ci}(x)\,\sin x- \operatorname{si}(x)\,\cos x =\int_{0}^{\infty}\frac{\sin t}{x+t}\,dt.
$$
Notice that $F\in \mathcal{S}_{2}\setminus \mathcal{S}_{1}$.

We note that when $\varphi$ is an even function on $\mathbb R$ satisfying the assumptions in Lemma \ref{lemma:continuous}, $\alpha>0$ and $0\leq\beta\leq T$ then 
\begin{align}
&\int_{0}^{\infty}\varphi(\alpha\,t+\beta)\,e^{-xt}\,dt=\nonumber \\
&\frac{1}{x}\left(\varphi(\beta)-\int_{0}^{\beta}\varphi^{\prime}(t)\,e^{-\frac{x(t-\beta)}{\alpha}}\,dt
+\,\,\frac{e^{\frac{\beta x}{\alpha}}}{1-e^{-\frac{Tx}{\alpha}}}
\,\,  \int_{0}^{T}\varphi^{\prime}(t)e^{\frac{-xt}{\alpha}}\,dt\right)\label{eq:+beta}
\end{align}
and
\begin{align}
&\int_{0}^{\infty}\varphi(\alpha\,t-\beta)\,e^{-xt}\,dt=\nonumber\\
&\frac{1}{x}\left(\varphi(\beta)-\int_{0}^{\beta}\varphi^{\prime}(t)\,e^{\frac{x(t-\beta)}{\alpha}}\,dt
+\,\,\frac{e^{-\frac{\beta x}{\alpha}}}{1-e^{-\frac{Tx}{\alpha}}}
\,\,  \int_{0}^{T}\varphi^{\prime}(t)e^{\frac{-xt}{\alpha}}\,dt\right).\label{eq:-beta}
\end{align}
The relations \eqref{eq:+beta} and \eqref{eq:-beta} give
\begin{align*}
&\frac{1}{2}\int_{0}^{\infty}\left(\varphi(\alpha\,t+\beta)+\varphi(\alpha\,t-\beta)\right)\,e^{-xt}\,dt=\\
&\frac{1}{x}\left(\varphi(\beta)-\int_{0}^{\beta}\varphi^{\prime}(t)\,\cosh\left(\frac{x}{\alpha}(t-\beta)\right)\,dt+\,\,\frac{\cosh\left(\frac{\beta}{\alpha}\,x\right)}{1-e^{\frac{-Tx}{\alpha}}}
\,\,  \int_{0}^{T}\varphi^{\prime}(t)e^{\frac{-xt}{\alpha}}\,dt\right).
\end{align*}
This proves that the function
\begin{align}
\lefteqn{\sigma(x)=}\label{eq:sigma-cont}\\
&\varphi(\beta)-\int_{0}^{\beta}\varphi^{\prime}(t)\,\cosh\left(\frac{x}{\alpha}(t-\beta)\right)\,dt+\,\,\frac{\cosh\left(\frac{\beta}{\alpha}\,x\right)}{1-e^{\frac{-Tx}{\alpha}}}
\,\,  \int_{0}^{T}\varphi^{\prime}(t)e^{\frac{-xt}{\alpha}}\,dt\nonumber
\end{align}
is infinitely divisible. Similarly, letting $A_{\varphi}=\max_{0\leq t\leq T}\{\varphi(t)\}$ and using again \eqref{eq:+beta} and \eqref{eq:-beta}, 
we obtain the infinite divisibility of 
\begin{align}
\lefteqn{\tau(x)=}\label{eq:cont}\\
&A_{\varphi}\pm 2\left(\int_{0}^{\beta}\varphi^{\prime}(t)\,\sinh\left(\frac{x}{\alpha}(t-\beta)\right)\,dt+\,\,\frac{\sinh\left(\frac{\beta}{\alpha}\,x\right)}{1-e^{\frac{-Tx}{\alpha}}}
\,\,  \int_{0}^{T}\varphi^{\prime}(t)e^{\frac{-xt}{\alpha}}\,dt\right).\nonumber
\end{align}

\begin{cor}
\label{cor:sigma-cont}
The function
$$
\sigma(x)=\nu-\frac{\sinh(\nu x)}{x}+\cosh(\nu x)\,\frac{\tanh (x)}{x}
$$
is infinitely divisible when $0\leq \nu\leq 1$. Its Laplace transform is in $\mathcal S_2$ and is given by
$$
L(\sigma)(x)=\frac{\nu}{x}+\log \left(\frac{\Gamma\left(\frac{x-\nu}{4}+1\right)\Gamma\left(\frac{x+\nu}{4}\right)}{\Gamma\left(\frac{x-\nu}{4}+\frac{1}{2}\right)\Gamma\left(\frac{x+\nu}{4}+\frac{1}{2}\right)}\right).
$$
\end{cor}
{\it Proof.} Applying \eqref{eq:sigma-cont} with $\varphi(t)=|t|$ for $-\pi\leq t\leq \pi$ and $T=2\pi$, $\alpha=\pi/2$, $\beta=\nu \pi/2$ 
yields the infinite divisibility of $\sigma$. To compute its Laplace transform
we first rewrite $\sigma$ as
$$
\sigma(x)=\nu+\frac{\sinh((1-\nu)x)}{x\,\cosh (x)}.
$$
It follows from Corollary \ref{cor:+-} (and Remark \ref{rem:psi-beta}) that 
  \begin{align*}
  \lefteqn{\int_{0}^{\infty}e^{-xt}\,\frac{\sinh((1-\nu) t)}{\cosh (t)}\,{\rm d}t}\\
  &=\frac{1}{2}\,\left(\beta\Big(\frac{x+\nu}{2}\Big)-\beta\Big(\frac{x-\nu}{2}+1\Big)\right)\\
  &=\frac{1}{4}\left(\psi\left(\frac{x+\nu}{4}+\frac{1}{2}\right)-\psi\left(\frac{x+\nu}{4}\right)-\psi\left(\frac{x-\nu}{4}+1\right)+\psi\left(\frac{x-\nu}{4}+\frac{1}{2}\right)\right)\\
  &=-\left[\log\left(\frac{\Gamma\big(\frac{x-\nu}{4}+1\big)\,\Gamma\big(\frac{x+\nu}{4}\big)}{\Gamma\big(\frac{x-\nu}{4}+\frac{1}{2}\big)\,\Gamma\big(\frac{x+\nu}{4}+\frac{1}{2}\big)}\right)\right]^{\prime}.
  \end{align*}
We certainly also have   
  \begin{align*}
  \int_{0}^{\infty}e^{-xt}\,\frac{\sinh((1-\nu) t)}{\cosh (t)}\,{\rm d}t&=-\left[\int_{0}^{\infty}e^{-xt}\,\frac{\sinh((1-\nu) t)}{t\,\cosh (t)}\,{\rm d}t \right]^{\prime}.
 \end{align*}
  Since the last two expressions in square brackets tend to zero as $x\to \infty$ the desired result follows.  \hfill $\square$

Of course many other examples can be obtained as applications of \eqref{eq:sigma-cont} and \eqref{eq:cont}. Let us give the analogue of Corollary \ref{cor:+-}.
The proof follows the same lines as the proof of Corollary \ref{cor:sigma-cont}.
\begin{cor}
The functions
$$
\tau(x)=1\pm\frac{1}{x}\left\{1-\frac{\cosh((1-\nu) x)}{\cosh(x)}\right\}
$$
are infinitely divisible when $0\leq \nu\leq 1$. The Laplace transforms belong to $\mathcal S_2$ and are given by
$$
L(\tau)(x)=\frac{1}{x}\pm\log \left(\frac{4\Gamma\left(\frac{x-\nu}{4}+1\right)\Gamma\left(\frac{x+\nu}{4}+\frac{1}{2}\right)}
{x\Gamma\left(\frac{x-\nu}{4}+\frac{1}{2}\right)\Gamma\left(\frac{x+\nu}{4}\right)}\right).
$$
\end{cor}

\section{Ratios of Gamma functions}
\label{sec:gamma} 
The Laplace transform in Corollary \ref{cor:sigma-cont} involves a ratio of Gamma functions. Such ratios have been intensively studied. In one of the first papers \cite{bustoz-ismail} the ratio 
 $$
 \frac{\Gamma(x)\Gamma(x+a+b)}{\Gamma(x+a)\Gamma(x+b)}
 $$
 was shown to be (logarithmically) completely monotonic when $a,b\geq 0$. In the next proposition we strengthen this result. 
\begin{prop} 
\label{prop:gamma-a-b}
Let $a,b> 0$. 
 The function  
 $$
 \log \frac{\Gamma(x)\Gamma(x+a+b)}{\Gamma(x+a)\Gamma(x+b)}
 $$
 belongs to $\mathcal S_2$  and is thus logarithmically completely monotonic. It has the representation
 $$
 \log \frac{\Gamma(x)\Gamma(x+a+b)}{\Gamma(x+a)\Gamma(x+b)}=\int_0^{\infty}\frac{g(t)}{(t+x)^2}\, dt,
 $$
 where
 $$
 g(t)=\sum_{k=0}^{[t]}\chi_{(k,\infty)}(t)\,(\chi_{(0,a)}*\chi_{(0,b)})(t-k).
 $$
\end{prop}

\begin{rem}
It follows that for $0<a<b$,
 $$
 \chi_{(0,a)}*\chi_{(0,b)}(t)=\left\{\begin{array}{ll}
              t,&0<t<a,\\
              a,&a<t<b,\\
              a+b-t,&b<t<a+b\\
              0,&t>a+b.
             \end{array}\right.
 $$           

\end{rem}
{\it Proof of Proposition \ref{prop:gamma-a-b}.}  The relation
\begin{equation}
\label{eq:bustoz-ismail}
\int_0^{\infty}e^{-xt}\frac{(1-e^{-at})(1-e^{-bt})}{t(1-e^{-t})}\, dt=\log \frac{\Gamma(x)\Gamma(x+a+b)}{\Gamma(x+a)\Gamma(x+b)}
\end{equation}
is implicitly contained in \cite{bustoz-ismail}, but for the reader's convenience we give the argument behind it.
By using the integral representation
$$
\psi(x)=-\gamma+\int_0^{\infty}\frac{e^{-t}-e^{-xt}}{1-e^{-t}}\, dt
$$
we immediately get
$$
-\int_0^{\infty}e^{-xt}\frac{(1-e^{-at})(1-e^{-bt})}{1-e^{-t}}\, dt=\left(\log \frac{\Gamma(x)\Gamma(x+a+b)}{\Gamma(x+a)\Gamma(x+b)}\right)'.
$$
The relation \eqref{eq:bustoz-ismail} now follows as in the proof of Corollary \ref{cor:sigma-cont}. From the elementary relations, letting $\mu=\sum_{k=0}^{\infty}\epsilon_k$,
$$
\frac{1-e^{-at}}{t}=L(\chi_{(0,a)})(t),\qquad \frac{1}{1-e^{-t}}=L(\mu)(t),
$$
we infer that the function 
$$
\frac{(1-e^{-at})(1-e^{-bt})}{t^2(1-e^{-t})}=L(\chi_{(0,a)}*\chi_{(0,b)}*\mu)(t)
$$
is completely monotonic. Furthermore, a computation shows that
$$
(\chi_{(0,a)}*\chi_{(0,b)}*\mu)(t)=\sum_{k=0}^{[t]}\chi_{(k,\infty)}(t)\,(\chi_{(0,a)}*\chi_{(0,b)})(t-k)=g(t).
$$
This completes the proof.\hfill $\square$

\begin{rem}
 In general the function in Proposition \ref{prop:gamma-a-b} does not belong to $\mathcal S_1$. However, when one of the parameters is an
 integer, say $b=n$, then it does. In that situation
 $$
 \frac{(1-e^{-at})(1-e^{-nt})}{t(1-e^{-t})}=\frac{1-e^{-at}}{t}\, \sum_{k=0}^{n-1}e^{-kt}
 $$
 is completely monotonic and 
 $$
 \log \frac{\Gamma(x)\Gamma(x+a+n)}{\Gamma(x+a)\Gamma(x+n)}=\log \frac{x+a}{x}+\cdots +\log \frac{x+a+n-1}{x+n-1}
 $$
 belongs to $\mathcal S_1$. This can of course also be seen directly as it is a sum of Stieltjes functions.  
\end{rem}

Let us mention an extension of Proposition \ref{prop:gamma-a-b} to entire functions of genus 1 having only negative zeros.

\begin{prop}
\label{prop:f-a-b}
 Let $f$ be an entire function of genus 1 having only negative zeros and such that $f(0)>0$.  The function
 $$
 \log \frac{f(x+a)f(x+b)}{f(x)f(x+c)}
 $$
 belongs to $\mathcal S_2$ when $a>0$, $b> 0$ and $c= a+b$. Its representing measure has density
 $\chi_{(0,a)}*\chi_{(0,b)}*\mu$ w.r.t.\ Lebesgue measure, where $\mu=\sum_{n=1}^{\infty}\epsilon_{a_n}$.
 
\end{prop}
{\it Proof.} Denoting the negative zeros by $\{-a_n\}$ with $0<a_1\leq a_2\leq \cdots$ we have (see \cite{pedersen})
$$
-(\log f)''(x)=\int_0^{\infty}e^{-xs}sh(s)\, ds,
$$
where $h(s)=\sum_{n=1}^{\infty}e^{-a_ns}$.
For $\alpha>0$ this yields
$$
-\left((\log f)'(x+\alpha)-(\log f)'(x)\right)=\int_0^{\infty}e^{-xs}(1-e^{-\alpha s})h(s)\, ds.
$$
Therefore we obtain
$$
-\left(\log \frac{f(x+a)f(x+b)}{f(x)f(x+c)}\right)'=\int_0^{\infty}e^{-xs}r(s)h(s)\, ds,
$$
with 
\begin{align*}
r(s)&=1-e^{-as}-e^{-bs}+e^{-cs}=(1-e^{-as})(1-e^{-bs}).
\end{align*}
Since $r(s)\sim abs^2$ as $s\to 0$, $h(s)r(s)/s$ is integrable at $0$ and we thus get
$$
\log \frac{f(x+a)f(x+b)}{f(x)f(x+c)}=\kappa+\int_0^{\infty}e^{-xs}sh(s)\frac{r(s)}{s^2}\, ds.
$$
It is easy to show that the constant $\kappa$ must be equal to $0$. (See also Remark \ref{rem:kappa}.) Furthermore, $r(s)/s^2$ is 
completely monotonic and this completes the proof. (The representing measure is found as in the proof of Proposition \ref{prop:gamma-a-b}.) \hfill $\square$
\begin{rem}
 In the situation where $a> 0$, $b> 0$ and $b\leq c\leq a+b$ the proof above shows that 
 $$
-\left(\log \frac{f(x+a)f(x+b)}{f(x)f(x+c)}\right)'=\int_0^{\infty}e^{-xs}s^2h(s)\frac{r(s)}{s^2}\, ds,
$$
with 
\begin{align*}
r(s)&
=(1-e^{-as})(1-e^{-bs})+e^{-cs}(1-e^{-(a+b-c)s}).
\end{align*}
belongs to $\mathcal S_3$. The function 
$$
\log \frac{f(x+a)f(x+b)}{f(x)f(x+c)}
$$
does not necessarily belong to $\mathcal S_2$.

When $a> 0$, $b> 0$ and $0\leq c\leq b$ we have 
$$
\log \frac{f(x+a)f(x+b)}{f(x)f(x+c)}=\log \frac{f(x+a)}{f(x)}+\log \frac{f(x+b)}{f(x+c)}
$$
is the negative of a Pick function (see Proposition \ref{prop:entire-again}). Since it tends to $-\infty$ as $x\to \infty$ it is not in any of the Stieltjes classes.
\end{rem}

Let us now investigate ratios of Gamma functions, leaving out the logarithm.

 \begin{prop} 
 \label{prop:gamma-N}
 For $N\in \{0,1,\ldots\}$ and $0<s<1$ the function  
 $$
 \frac{\Gamma(x)}{\Gamma(x+N+s)}
 $$
 belongs to $\mathcal S_{N+1}$.
 \end{prop}
 {\it Proof.}  For $s>0$ we have
 \begin{equation}
 \label{eq:gammaxs}
 \frac{\Gamma(x)}{\Gamma(x+s)}=\frac{1}{\Gamma(s)}\int_0^{\infty}e^{-xt}(1-e^{-t})^{s-1}\, dt,\quad x>0.
 \end{equation}
 Using the fact that $1/(1-e^{-t})$ is logarithmically completely monotonic it follows that $(1-e^{-t})^{s-1}$ is completely monotonic when $0<s<1$. Thus $\Gamma(x)/\Gamma(x+s)$ is in
 $\mathcal S_1$ for $s\in (0,1)$.
 
 The functional equation for the Gamma function entails
 $$
 \frac{\Gamma(x)}{\Gamma(x+N+s)}=\frac{1}{x}\frac{1}{x+1}\cdots\frac{1}{x+N-1}\frac{\Gamma(x+N)}{\Gamma(x+N+s)}.
 $$
 Each term in this product belongs to $\mathcal S_1$, and hence the product belongs to $\mathcal S_{N+1}$.\hfill $\square$
 \begin{rem}
  When $N=1$ the proposition states that $\Gamma(x)/\Gamma(x+s+1)\in\mathcal S_2$ for $0<s<1$. (And it does not belong to $\mathcal S_1$.) Since 
  $$
  L\left(\sum_{k=0}^{\infty}\frac{(1-s)_k}{k!}\chi_{(k,k+1)}\right)(t)=\frac{(1-e^{-t})^s}{t}
  $$
  it has the integral representation
  $$
  \frac{\Gamma(x)}{\Gamma(x+s+1)}=\frac{1}{\Gamma(s+1)}\int_0^{\infty}\frac{1}{(x+t)^2}\sum_{k=0}^{\infty}\frac{(1-s)_k}{k!}\chi_{(k,k+1)}(t)\, dt,\quad x>0.
  $$
 \end{rem}

 \begin{rem} 
 We cannot expect that $\Gamma(x)/\Gamma(x+s)$ belongs to $\mathcal S_s$ for all $s\in (0,1)$. Indeed, if this were the case then \eqref{eq:gammaxs} would imply that
 the function 
 $$
 \frac{t}{1-e^{-t}}
 $$
 would be logarithmically completely monotonic. This is however not the case. (It is not even decreasing!) 
 \end{rem}
We notice that there is a different way of proving Proposition \ref{prop:gamma-N} using the so-called Pick (or Nevanlinna) functions. A Pick function is by definition a holomorphic function 
with non-negative imaginary part in the upper half-plane. Any Pick function $p$ admits an integral representation of the form
$$
p(z)=\alpha z+\beta+\int_{-\infty}^{\infty}\left(\frac{1}{t-z}-\frac{t}{t^2+1}\right)\, d\mu(t),
$$
where $\alpha\geq 0$, $b\in \mathbb R$ and $\mu$ satisfies $\int_{-\infty}^{\infty}d\mu(t)/(1+t^2)<\infty$.
Furthermore, $\alpha=\lim_{y\to \infty}p(iy)/(iy)$, $\beta=\Re p(i)$ and $\mu=1/\pi \lim_{\epsilon\to 0_+}\Im p(t+i\epsilon)dt$ vaguely. Let $A\subset \mathbb R$ be closed.
The function $p$ has a holomorphic extension to $\mathbb C\setminus A$ satisfying $p(\overline{z})=\overline{p(z)}$ if and only if the support of $\mu$ is contained in $A$. 
Here is a simple observation.
 \begin{prop}
   \label{prop:entire-again}
 Let $f$ be a real entire function of genus $1$ with only negative zeros and with $f(0)>0$. For $s>0$ the function 
 $$
 h(z)=(\log f)(z)-(\log f)(z+s),\quad z\in \mathbb C\setminus (-\infty,0]
 $$
 is a Pick function. Its representing measure is supported on $(-\infty,0]$  and $\lim_{y\to \infty}h(iy)/(iy)=0$.
 \end{prop}

\noindent
 {\it Proof.} From Hadamard's factorization theorem it follows that
$$
\log f(z)=az+b+\sum_{n=1}^{\infty}\left(\Log\left(1+\frac{z}{a_n}\right) -\frac{z}{a_n}\right),
$$
where $\Log$ on the right hand side denotes the principal logarithm, and  $\{-a_n\}$ are the zeros of $f$ arranged so that $0<a_1\le a_2\le\ldots$. Hence
$$
h(z)=-as+\sum_{n=1}^{\infty}\left(\Log \left(1+\frac{z}{a_n}\right)-\Log \left(1+\frac{z+s}{a_n}\right)+\frac{s}{a_n}\right)
$$
and this gives 
$$
\Im h(z)=\sum_{n=1}^{\infty}\left(\Arg \left(1+\frac{z}{a_n}\right)-\Arg \left(1+\frac{z+s}{a_n}\right)\right).
$$
For $\Im z>0$ each summand is positive. Therefore $h$ is a Pick function. To prove the second assertion it is convenient to write the Hadamard representation in the form
$$
\log f(z)=\alpha z+\beta+\int_0^{\infty}( \Log (1+z/t)-z/t)\, dn(t),
$$
where $n$ denotes the zero counting function $n(t)=\#\{n\geq 1\,|\, a_n\leq t\}$. From this we obtain by partial integration
$$
\log f(z)=\alpha z+\beta-z^2\int_0^{\infty}\frac{n(t)}{t^2}\frac{dt}{t+z}
$$
so that 
\begin{align*}
 h(iy)&=-\alpha s+sy^2\int_0^{\infty}\frac{n(t)}{t^2}\frac{dt}{(t+iy)(t+s+iy)}\\
 &\phantom{=}+2iys\int_0^{\infty}\frac{n(t)}{t^2}\frac{dt}{t+s+iy}+s^2\int_0^{\infty}\frac{n(t)}{t^2}\frac{dt}{t+s+iy}.
\end{align*}
It now follows that $\lim_{y\to \infty}h(iy)/(iy)=0$ by Lebesgue's theorem on dominated convergence. This completes the proof.\hfill $\square$

It is evident that the derivative of a Pick function, whose representing measure is supported on $(-\infty, 0]$, belongs to $\mathcal S_2$. 
This gives the proof of 
the next observation.
\begin{cor}
 \label{prop:entire}
 Let $f$ be a real entire function of genus $1$ with only negative zeros and with $f(0)>0$. For $s>0$ the function 
 $$
 h(x)=(\log f)'(x)-(\log f)'(x+s),\quad x>0
 $$
 belongs to  $\mathcal S_2$.
\end{cor}

\begin{rem}
\label{rem:kappa}
For a Pick function $p$ with integral representation 
$$
p(z)=\beta+\int_{-\infty}^{0}\left(\frac{1}{t-z}-\frac{t}{t^2+1}\right)\, d\mu(t)
$$
we obtain
\begin{align*}
p(x)-p(x+s)&=\int_{-\infty}^{0}\left(\frac{1}{t-x}-\frac{1}{t-x-s}\right)\, d\mu(t)\\
&=-s\int_0^{\infty}\frac{d\mu(-t)}{(t+x)(t+x+s)}.
\end{align*}
By dominated convergence it follows that
$\lim_{x\to \infty}(p(x)-p(x+s))=0$.

The function in Proposition \ref{prop:f-a-b} can be written as $p(x+b)-p(x)$, where $p(x)=\log f(x)-\log f(x+a)$. It is thus a difference of two Pick functions of the form in Proposition \ref{prop:entire-again}. Hence the constant $\kappa$ in the 
proof of Proposition \ref{prop:f-a-b} equals zero.
\end{rem}

Let us finally state the integral representation of the Pick function in Proposition \ref{prop:entire-again} 
in the case where $f=1/\Gamma$ is the reciprocal of the Gamma function, and $s\in(0,1)$. (Since $f$ has a zero at the origin the proof of that proposition needs to be modified 
slightly.)
$$
\log \Gamma(z+s)-\log \Gamma(z)=\int_{-\infty}^0\left(\frac{1}{t-z}+\frac{t}{t^2+1}\right)\, d\tau(t)+c_s
$$
where $c_s\in \mathbb R$ and $\tau$ is the measure
$$
\tau=\sum_{k=0}^{\infty}dt|_{(-k-s,-k)}.
$$
It now follows that for $\Im z>0$, the imaginary part of $\log \Gamma(z+s)-\log \Gamma(z)$ is less than 
$$
\int_{-\infty}^0\frac{y}{(t-x)^2+y^2}\, dt,
$$ 
which is less than $\pi$. This on the other hand actually shows that $\Gamma(z+s)/\Gamma(z)$ is also a Pick function and from here it follows (using the complex characterization
of Stieltjes functions) that the reciprocal $\Gamma(z)/\Gamma(z+s)$ is in fact in $\mathcal S_1$.

\begin{rem}
 The Pick property of the ratio of two Gamma functions is a special case of a result due to Krein (see e.g.\ \cite[Section 27.2]{levin}) concerning meromorphic functions with interlacing zeros and poles
 mapping the upper half-plane into itself.
\end{rem}

\section{The double Gamma function of Barnes}
\label{sec:barnes}
The results in the previous sections have implications for the remainders in asymptotic expansions of the double Gamma function of Barnes. 

Asymptotic expansions of the (Euler) Gamma
function is a classical and much investigated area. Stirling's formula and series is usually expressed in terms of $\log \Gamma$. Monotonicity properties of the remainders 
in these expansions have also been studied. The well-known relation (see \cite[p.~57]{remmert})
$$
\log \Gamma(z)-\left((z-\tfrac{1}{2})\log z-z+\tfrac{1}{2}\log (2\pi)\right)=\int_0^{\infty}\frac{Q(t)}{(z+t)^2}\, dt, 
$$
for $z\in \mathbb C\setminus (-\infty, 0]$ and where $Q(t)=\tfrac{1}{2}(t-[t]-(t-[t])^2)$, can be restated as the fact that the remainder in the approximation of $\log \Gamma$ by the elementary function in the 
parentheses belongs to $\mathcal S_2$.

Barnes (see \cite {barnes} and \cite{barnes1}) introduced the so-called double gamma function (and multiple gamma functions) more than 100 years ago. 
For a modern introduction we refer to \cite{ruijsenaars}. There, an asymptotic expansion of in particular the logarithm of the double gamma function was obtained. It reads
\begin{align*}
\log \Gamma_2(w)&= -\frac{B_{2,2}(w)}{2}\log w +
\frac{3}{4}B_{2,0}(0)w^2+B_{2,1}(0)w\\
&\phantom{=} +\sum_{k=3}^{m}\frac{(-1)^k}{k!}(k-3)!B_{2,k}(0)w^{2-k} + R_{2,m}(w),
\end{align*}
where $R_{2,m}$ is a certain remainder. Here, $\{B_{2.k}(x)\}$ is the sequence of double Bernoulli polynomials defined via
$$
\frac{t^2e^{xt}}{(e^t-1)^2}=\sum_{k=0}^{\infty}B_{2,k}(x)\frac{t^k}{k!}.
$$
For $n\geq 1$, the even indexed remainder $R_{2,2n}$ was shown (see \cite{pedersen1, kp1}) to have the representation
$$
R_{2,2n}(w)= \int_{0}^{\infty}e^{-wt}t^{2n}p_n(t)\, dt,
$$
where 
\begin{equation}
\label{eq:p_n}
\begin{array}{rl}
p_n(t)=\displaystyle{\frac{1}{t^2}\sum_{k=1}^{\infty}(2\pi k)^{1-2n}} &\left(  \displaystyle{\frac{4\pi k}{t^2+(2\pi k)^2}+\frac{8\pi kt}{(t^2+(2\pi k)^2)^2}+}\right. \\
& \left.  \displaystyle{\frac{(2n-1)}{2\pi k}\frac{2t}{t^2+(2\pi
    k)^2}}\right).
\end{array}
\end{equation}
\begin{prop}
 The function $p_n$ defined in \eqref{eq:p_n} is completely monotonic for any $n\geq 1$.
\end{prop}
{\it Proof.} It is enough to verify that 
$$
f_{k,n}(t)=\frac{1}{t^2}\left(\frac{4\pi k}{t^2+(2\pi k)^2}+\frac{8\pi kt}{(t^2+(2\pi k)^2)^2}+
\frac{(2n-1)}{2\pi k}\frac{2t}{t^2+(2\pi k)^2}\right)
$$
is completely monotonic, and by making the substitution $t=2\pi k x$, this amounts to showing that
$$
g_{k,n}(x)=\frac{1}{x^2(x^2+1)}+\frac{1}{\pi k}\frac{1}{x(x^2+1)^2}+\frac{2n-1}{2\pi k}\frac{1}{x(x^2+1)}
$$
is completely monotonic. Here, the first and third term is completely monotonic, but the second is not and we thus need to combine the terms and take cancellation into account.
Complete monotonicity of $g_{k,n}$ will follow from the complete monotonicity of 
$$
g(x)=\frac{1}{x^2(x^2+1)}+\frac{1}{\pi k}\frac{1}{x(x^2+1)^2},
$$
which we now turn to verify.

It is easy to show that 
$$
g(x)=L\left((1-\cos t)*1+\tfrac{1}{\pi k}(1-\cos t)*\sin t\right)(x)
$$
and that 
$$
(1-\cos t)*1+\tfrac{1}{\pi k}(1-\cos t)*\sin t=t-\sin t+\tfrac{1}{\pi k}(1-\cos t-\tfrac{1}{2}t\sin t).
$$
According to Lemma \ref{lemma:pos} this expression is positive for $t>0$. The proof is complete. \hfill $\square$
\begin{cor}
 The function  
 $$
 R_{2,2}(x)=\log \Gamma_2(x)+\left(\tfrac{1}{2}x^2-x+\tfrac{5}{12}\right)\log x -\tfrac{3}{4}x^2+x
 $$
 belongs to $\mathcal S_3$.
\end{cor}
\begin{lemma}
 \label{lemma:pos}
 For any $c\in [0,1]$ the function 
 $$h(t)=t-\sin t+c(1-\cos t-\tfrac{1}{2}t\sin t),\quad t>0
 $$
 is non-negative.
\end{lemma}
{\it Proof.} Standard computation shows $h$ is increasing on $(0,\pi)$ and for $t\geq \pi$ we use
$$
h(t)\geq t-\sin t- \tfrac{1}{2}ct\sin t\geq (1-\tfrac{1}{2}c)t-1>0,
$$
taking into account $0\leq c\leq 1$.\hfill $\square$

\section{Some generalized Stieltjes functions of special type}
\label{sec:cesaro}
In this section we consider sums of the form 
\begin{equation}
\label{eq:gen}
\sum_{n=0}^{\infty}\frac{a_n}{(x+n)^{\lambda}},
\end{equation}
where certain iterated partial sums of the real numbers $a_n$ are nonnegative. 

We remark that the series in \eqref{eq:gen}
is related to the so-called generalized Mathieu series, considered, for example, in \cite{pogany} and \cite{Qi}, where among
other things an integral representation was derived. 

For a real sequence $\{a_n\}_{n\geq 0}$ let 
$$
s^{(0)}_n=\sum_{j=0}^na_j
$$
and let for $k\geq 1$
\begin{equation}
\label{eq:sk}
s^{(k)}_n=\sum_{j=0}^ns^{(k-1)}_j.
\end{equation}
The sequence $\{s^{(0)}_n\}$ is the sequence of partial sums of the perhaps divergent series $\sum_{n=0}^{\infty}a_n$. A relation between the iterated partial 
sums is given in the next lemma. 
\begin{lemma}
\label{lemma:s}
For $k\geq 0$,
 $$
 (1-x)^{k+1}\sum_{n=0}^Ns^{(k)}_nx^n=\sum_{n=0}^{N}a_nx^n-x^{N+1}\sum_{j=0}^{k}s^{(j)}_N(1-x)^{j}.
 $$
\end{lemma}
{\it Proof.} The relations
\begin{align*}
 (1-x)\sum_{n=0}^Ns^{(0)}_nx^n&=\sum_{n=0}^{N}a_nx^n-s^{(0)}_Nx^{N+1}, \\
 (1-x)\sum_{n=0}^Ns^{(k)}_nx^n&=\sum_{n=0}^{N}s^{(k-1)}_nx^n-s^{(k)}_Nx^{N+1},\quad k\geq 1
\end{align*}
 follow easily using \eqref{eq:sk}. From the latter of these we infer (for $k\geq 1$)
 \begin{align*}
  (1-x)^{k+1}\sum_{n=0}^Ns^{(k)}_nx^n&=(1-x)^{k}\sum_{n=0}^Ns^{(k-1)}_nx^n-s^{(k)}_N(1-x)^kx^{N+1},
 \end{align*}
and combining with the former relation this completes the proof.\hfill$\square$

The proposition below gives conditions under which the expression in \eqref{eq:gen} defines a generalized Stieltjes functions of positive order. 
\begin{prop}
\label{prop:cesaro}
 Suppose that for some $k\geq 0$ we have 
 \begin{enumerate}[(i)]
 \item $s^{(k)}_n\geq 0$ for all $n\geq 0$,
 \item $s^{(k)}_n/n^{\lambda}\to 0$ for $n\to \infty$,
 \item $\sum_{n=1}^{\infty}s^{(k)}_n/n^{1+\lambda}<\infty$.
 \end{enumerate}
 Then the sum
$$
f(x)=\sum_{n=0}^{\infty}\frac{a_n}{(x+n)^{\lambda}},\quad x>0
$$
converges and $f$ belongs to $\mathcal S_{\lambda+k+1}$. It has the representation
$$
f(x)=\int_0^{\infty}\frac{d\mu(t)}{(x+t)^{\lambda+k+1}}=\frac{1}{\Gamma(\lambda)}\int_{0}^{\infty} e^{-xt}\, t^{\lambda+k} \kappa(t)\,dt,
$$
where 
$$
d\mu(t)=\frac{(\lambda)_{k+1}}{k!}\,\sum_{n=0}^{\infty}a_{n}\,\chi_{(n,\infty)}(t)(t-n)^{k}dt
$$ 
is a positive measure, and
$$
\kappa(t)=\frac{1}{t^{k+1}}\sum_{n=0}^{\infty} a_{n}\,e^{-nt}
$$
is completely monotonic.
\end{prop}
\begin{rem}
If condition (ii) in Proposition \ref{prop:cesaro} holds then also  $s^{(k-1)}_n/n^{\lambda}=s^{(k)}_n/n^{\lambda}-s^{(k)}_{n-1}/n^{\lambda}$ tends 
to zero as $n$ tends to infinity. For any $j$ with $0\leq j\leq k$ we thus have  
$$
s^{(j)}_n/n^{\lambda}\to 0 \quad \text{for}\quad n\to \infty.
$$
\end{rem}
\begin{rem}
The measure $\mu$ in Proposition \ref{prop:cesaro} has the density given by 
$$
\frac{1}{k!}\,\sum_{j=0}^{m}a_{j}(t-j)^k, 
$$
for $t\in(m,m+1)$ and $m=0,1,\ldots$.
\end{rem}

\noindent
{\it Proof of Proposition \ref{prop:cesaro}.} The power series
$$
\sum_{n=0}^{\infty}s^{(k)}_nz^n
$$
has radius of convergence at least $1$ (because of (ii)). The power series $\sum_{n=0}^{\infty}a_nz^n$ has the same radius of convergence and we have (by Lemma \ref{lemma:s})
\begin{equation}
\label{eq:cb}
\frac{1}{(1-z)^{k+1}}\sum_{n=0}^{\infty}a_{n}z^{n}=\sum_{n=0}^{\infty}s^{(k)}_{n}z^{n},\quad |z|<1.
\end{equation}
Let now 
$$
g(t)=\sum_{n=0}^{\infty}s^{(k)}_n\chi_{(n,n+1)}(t).
$$
The Laplace transform of $g$ is defined and equals
$$
L(g)(x)=\frac{1-e^{-x}}{x}\sum_{n=0}^{\infty}s^{(k)}_ne^{-nx}.
$$
This gives
\begin{equation}
\label{eq:g*}
L(g*(\chi_{(0,1)})^{*k})(x)=\left(\frac{1-e^{-x}}{x}\right)^{k+1}\sum_{n=0}^{\infty}s^{(k)}_ne^{-nx}. 
\end{equation}
The function $t^{\lambda+k}L(g*(\chi_{(0,1)})^{*k})(t)$ is by \eqref{eq:cb} continuous on $(0,\infty)$. It is polynomially bounded for $t\to \infty$. 
Indeed,
$$
t^{\lambda+k}L(g*(\chi_{(0,1)})^{*k})(t)=t^{\lambda-1}(1-e^{-t})^{k+1}\sum_{n=0}^{\infty}s^{(k)}_{n}e^{-nt},
$$
and the infinite sum is a deceasing function of $t$. In order to show (also using that $1-e^{-t}\sim t$ for $t\to 0$) that the Laplace transform of the function in \eqref{eq:g*} is defined it is 
thus enough to show that 
$$
t^{\lambda+k}\sum_{n=0}^{\infty}s^{(k)}_{n}e^{-nt}
$$
is integrable on e.g.\ $(0,\log 2)$: 
\begin{align*}
\int_0^{\log 2}t^{\lambda+k}\sum_{n=1}^{\infty}s^{(k)}_{n}e^{-nt}\,dt&=\int_{1/2}^1(-\log u)^{\lambda+k}\sum_{n=1}^{\infty}s^{(k)}_{n}u^{n-1}\, du\\
&=\sum_{n=1}^{\infty}s^{(k)}_{n}\int_{1/2}^1(-\log u)^{\lambda+k}u^{n-1}\, du\\
&\leq \sum_{n=1}^{\infty}\left(2\log 2\right)^{\lambda+k}s^{(k)}_{n}\int_{1/2}^1(1-u)^{\lambda+k}u^{n-1}\, du\\
&\leq \left(2\log 2\right)^{\lambda+k}\sum_{n=1}^{\infty}s^{(k)}_{n}B(\lambda+k+1,n).
\end{align*}
Since $B(\lambda+k+1,n)\sim \Gamma(\lambda+k+1)n^{-\lambda-k-1}$ as $n$ tends to infinity it follows that the sum converges, and thus that the Laplace transform 
$$
\int_0^{\infty}e^{-xt}t^{\lambda+k}L(g*(\chi_{(0,1)})^{*k})(t)\, dt
$$
is defined and belongs to $\mathcal S_{\lambda+k+1}$. Now, by monotone convergence, 
\begin{align*}
\lefteqn{\int_0^{\infty}e^{-xt}t^{\lambda+k}L(g*(\chi_{(0,1)})^{*k})(t)\, dt}\\
&=\int_0^{\infty}e^{-xt}t^{\lambda-1}(1-e^{-t})^{k+1}\sum_{n=0}^{\infty}s^{(k)}_{n}e^{-nt}\, dt\\
&=\lim_{N\to \infty}\int_0^{\infty}e^{-xt}t^{\lambda-1}(1-e^{-t})^{k+1}\sum_{n=0}^{N}s^{(k)}_{n}e^{-nt}\, dt.
\end{align*}
Here, by Lemma \ref{lemma:s},
\begin{align*}
 \lefteqn{\int_0^{\infty}e^{-xt}t^{\lambda-1}(1-e^{-t})^{k+1}\sum_{n=0}^{N}s^{(k)}_{n}e^{-nt}\, dt}\\
 &=\sum_{n=0}^{N}a_{n}\int_0^{\infty}e^{-xt}t^{\lambda-1}e^{-nt}\, dt-
 \sum_{j=0}^ks^{(j)}_N\int_0^{\infty}e^{-xt}t^{\lambda-1}e^{-(N+1)t}(1-e^{-t})^{j}\,dt\\
 &=\Gamma(\lambda)\sum_{n=0}^{N}\frac{a_{n}}{(n+x)^{\lambda}}-
 \sum_{j=0}^k\frac{s^{(j)}_N}{(x+N+1)^{\lambda}}\int_0^{\infty}e^{-s}s^{\lambda-1}(1-e^{-s/(x+N+1)})^{j}\,ds.
 \end{align*}
Consider now the $j$th term in the last sum. The factor $s^{(j)}_N/(x+N+1)^{\lambda}$ tends to $0$  as $N\to \infty$, and the integral factor remains bounded as $N\to \infty$.
We conclude that  
$$
\Gamma(\lambda)\sum_{n=0}^{\infty}\frac{a_{n}}{(n+x)^{\lambda}}=\int_0^{\infty}e^{-xt}t^{\lambda+k}L(g*(\chi_{(0,1)})^{*k})(t)$$
belongs to $\mathcal S_{\lambda+k+1}$. Furthermore, $\kappa$ defined by $\kappa=L(g*(\chi_{(0,1)})^{*k})$ satisfies
$$
\kappa(t)=\frac{1}{t^{k+1}}\sum_{n=0}^{\infty} a_{n}\,e^{-nt}.
$$
The measure $\mu$ representing $f$ as a generalized Stieltjes function of order $\lambda+k+1$ is given by
$$
d\mu(t)=g*(\chi_{(0,1)})^{*k}(t)dt.
$$
Finally, from \eqref{eq:g*} we obtain that
\begin{align*}
L(g*(\chi_{(0,1)})^{*k})(x)&=\lim_ {N\to \infty}\frac{1}{k!}\sum_{n=0}^Na_nL(t^k*\epsilon_n)(x)\\
&=\lim_ {N\to \infty}\frac{1}{k!}L\left(\sum_{n=0}^Na_n\chi_{(n,\infty)}(t)(t-n)^k\right)(x),
\end{align*}
so that 
$$
g*(\chi_{(0,1)})^{*k}(t)=\frac{1}{k!}\,\sum_{n=0}^{\infty}a_{n}\,\chi_{(n,\infty)}(t)(t-n)^{k}.
$$
This completes the proof.\hfill $\square$

Let us end the paper by considering two applications of Proposition \ref{prop:cesaro}.
\begin{ex}
 Let $\lambda>0$ and $0<a\leq 1$ and consider the function
$$
\beta_{a, \lambda}(x)=\frac{1}{\Gamma(\lambda)}\int_{0}^{\infty} e^{-xt}\,\frac{1}{(1+e^{-t})^a}\, t^{\lambda-1} \,dt=\sum_{n=0}^{\infty} (-1)^{n}\frac{(a)_{n}}{n!}\, \frac{1}{(x+n)^{\lambda}}.
$$
Since, for $0<a\leq 1$,
$$
s_n^{(0)}=\sum_{j=0}^{n}(-1)^{j}\frac{(a)_{j}}{j!}\geq 0,\quad n=0,1,\ldots
$$
and $\{s_n^{(0)}\}$ is also a bounded sequence, it follows that $\beta_{a, \lambda}(x)$ is in $\mathcal S_ {\lambda+1}$. In particular, Nielsen's beta function $\beta(x)=\beta_{1, 1}(x)$ is in $\mathcal S_ {2}$.
\end{ex}

\begin{ex}
Prym's function, see \cite[p.~27]{N}, is given by
$$
P(x)=\sum_{n=0}^\infty\frac{(-1)^n}{n!(x+n)},
$$
and for $x>0$ it enters in the decomposition $\Gamma(x)=P(x)+Q(x)$, with
$$
P(x)=\int_0^1 t^{x-1}e^{-t}\,dt,\quad Q(x)=\int_1^\infty t^{x-1}e^{-t}\,dt.
$$
(These functions are special cases of the functions $P_a$, $Q_a$ with $a=1$, see \cite[Kapitel II]{N}.)

The function $P$ is a special case of (18) with $\lambda=1$ and $a_n=(-1)^n/n!$.
It fits into Proposition \ref{prop:cesaro} with $k=0$ since
$$
s_n^{(0)}=\sum_{j=0}^n \frac{(-1)^j}{j!}\ge 0, \quad n=0,1,\ldots,
$$
and $\{s_n^{(0)}\}$ is also bounded. Prym's function thus belongs to $\mathcal S_2$, and the corresponding representations hold with 
$$
d\mu(t)=\left(\sum_{n<t} \frac{(-1)^n}{n!}\right) \,dt
$$
and
$$
\kappa(t)=\frac{1}{t}\sum_{n=0}^\infty\frac{(-1)^n}{n!}e^{-nt}=\frac{1}{t}\exp(-e^{-t}).
$$
Hence
$$
P(x)=\int_0^\infty e^{-xt}\exp(-e^{-t})\,dt =\int_0^1 u^{x-1}e^{-u}\,du,
$$
after the substitution $u=e^{-t}$.

In particular, the probability density on $(0,\infty)$ given by 
$$
d_a(t)=\frac{1}{P(a)}\exp(-at-e^{-t})
$$
is infinitely divisible. For $a=1$ where $P(1)=1-e^{-1}$ we obtain the infinite divisibility of the half-Gumbel density, cf.\ \cite[p.~412]{steutel-van-harn}.
\end{ex}

\noindent
Christian Berg\\
Department of Mathematical Sciences\\
University of Copenhagen \\
Universitetsparken 5\\
DK-2100, Denmark\\
{\em email}:\hspace{2mm}{\tt berg@math.ku.dk}

\vspace{0.5cm}

\noindent
Stamatis Koumandos\\
Department of Mathematics and Statistics\\
The University of Cyprus \\ 
P.\ O.\ Box 20537\\
1678 Nicosia, Cyprus\\
{\em email}:\hspace{2mm}{\tt skoumand@ucy.ac.cy}

\vspace{0.5cm}

\noindent
Henrik Laurberg Pedersen\\
Department of Mathematical Sciences\\
University of Copenhagen \\
Universitetsparken 5\\
DK-2100, Denmark\\
{\em email}:\hspace{2mm}{\tt henrikp@math.ku.dk}

\end{document}